# Optimal Design Measures under Asymmetric Errors, with Application to Binary Design Points

Mausumi Bose  and  Rahul Mukerjee
Applied Statistics Division         Indian Institute of Management Calcutta
Indian Statistical Institute        Joka, Diamond Harbour Road
203 B.T. Road, Kolkata 700 108, India   Kolkata 700 104, India

**Abstract**: We study the optimal design problem under second-order least squares estimation which is known to outperform ordinary least squares estimation when the error distribution is asymmetric. First, a general approximate theory is developed, taking due cognizance of the nonlinearity of the underlying information matrix in the design measure. This yields necessary and sufficient conditions that a *D*- or *A*-optimal design measure must satisfy. The results are then applied to find optimal design measures when the design points are binary. The issue of reducing the support size of the optimal design measure is also addressed.

## 1. Introduction

The existing literature on optimal designs focuses on maximizing or minimizing some scalar function of the covariance matrix of the ordinary least squares (OLS) estimator. In recent years, however, it has been well recognized (Wang and Leblanc, 2008) that the second-order least squares (SLS) estimator is more efficient than the OLS estimator when the error distribution is asymmetric. Another advantage of SLS estimation is that it is applicable even when the error distribution is not exactly known. Gao and Zhou (2014) pioneered research on the optimal design problem under SLS estimation. They formulated the *D*- and *A*-optimality criteria with reference to the asymptotic covariance matrix of SLS estimators of the parameters of interest and then explored optimal regression designs under these criteria. See their paper also for an excellent review of SLS estimation and its advantages over OLS estimation, together with further references.

   In the present work, we first develop a general approximate theory for the optimal design problem under SLS estimation and obtain necessary and sufficient conditions characterizing *D*- and *A*-optimal design measures. This is done in Section 2 taking due cognizance of the nonlinearity of the underlying information matrix in the design measure. While the resulting optimality conditions look more involved than their counterparts under OLS estimation, they have a broad spectrum of applicability allowing analytical as well as algorithmic derivation of optimal design measures in specific contexts. This is illustrated at length in Section 3 with reference to design spaces with binary design points as arise, for example, in spring balance weighing designs. The issue of reducing the support size of the optimal design measure is also addressed.

Our approximate theory results in Section 2 are given in the framework of a linear model and a finite design space for ease in notation and presentation, and also because this is what we need in the application that follows. Their message is, however, quite general. With heavier notation, the same arguments readily yield versions of these results characterizing locally optimal design measures in non-linear models under SLS estimation. Moreover, it is not hard to extend these results to continuous design spaces, with sums replaced by integrals as and when needed.

## 2. Approximate theory

Consider a linear model such that the expectation of any observation is of the form $x^T\theta$, where $\theta$ is a $q$-dimensional unknown parameter, $x$ is a $q\times 1$ design point which belongs to a design space $\Omega$ and the superscript $T$ denotes transpose. We assume as usual that the observational errors have a common variance and are uncorrelated, but allow them to have a possibly asymmetric distribution and study optimal design measures for $\theta$ under SLS estimation. Suppose the design space $\Omega$ is finite, consisting of design points $x_1,...,x_n$, each a $q\times 1$ vector. A design measure $p = (p_1,..,p_n)^T$ assigns masses $p_1,..,p_n\,(\geq 0)$ to $x_1,...,x_n$, respectively, where $p_1+..+p_n=1$. Let

$$G(p)=\Sigma_{i=1}^n p_i x_i x_i^T, \qquad g(p)=\Sigma_{i=1}^n p_i x_i, \qquad H(p)=G(p)-tg(p)g(p)^T, \qquad (1)$$

where $t=\mu_3^2/\{\mu_2(\mu_4-\mu_2^2)\}$, with $\mu_j$ representing the $j$th central moment of the error distribution. Clearly, $0\leq t<1$.

Following Gao and Zhou (2014), the asymptotic covariance matrix of the SLS estimator of $\theta$, under a design measure $p$, is proportional to $\{H(p)\}^{-1}$. Hence we study optimal design measures on the basis of criteria defined with respect to $H(p)$ which may be viewed as the information matrix for $\theta$ under SLS estimation, analogously to $G(p)$ under OLS estimation. Thus the $D$- and $A$-criteria now call for finding $p$ so as to maximize $\phi_D(p)$ and $\phi_A(p)$, respectively, where

$$\phi_D(p)=\log[\det\{H(p)\}], \quad \text{if } H(p) \text{ is nonsingular,}$$
$$= -\infty, \qquad \text{otherwise,} \qquad (2)$$

$$\phi_A(p)=-\operatorname{tr}\{H(p)\}^{-1}, \quad \text{if } H(p) \text{ is nonsingular,}$$
$$= -\infty, \qquad \text{otherwise.} \qquad (3)$$

In contrast to $G(p)$ arising under OLS estimation, the matrix $H(p)$ is nonlinear in $p$ for $t>0$. Lemmas 1-4 below show that still the basic arguments of approximate theory, that make use of directional derivatives, remain valid in our setup. Lemma 3 gives a theoretical characterization of $D$- and $A$-optimal design measures while Lemma 4 will be seen to be useful in their algorithmic determination. Lemmas 1 and 2 prepare the background for Lemmas 3 and 4. To avoid trivialities, here-after, we assume that for each $t$, there exists a design measure $p$ with nonsingular $H(p)$. For any such $p$ and $1 \leq i \leq n$, let

$$\psi_{Di}(p) = (1-t)x_i^T \{H(p)\}^{-1} x_i + t\{x_i - g(p)\}^T \{H(p)\}^{-1}\{x_i - g(p)\}, \qquad (4)$$

$$\psi_{Ai}(p) = (1-t)x_i^T \{H(p)\}^{-1}\{H(p)\}^{-1} x_i + t\{x_i - g(p)\}^T \{H(p)\}^{-1}\{H(p)\}^{-1}\{x_i - g(p)\}. \qquad (5)$$

**Lemma 1**. *Both $\phi_D(p)$ and $\phi_A(p)$ are concave in p.*

**Lemma 2**. *Let $\tilde{p} = (\tilde{p}_1, ..., \tilde{p}_n)^T$ and p be design measures such that $H(p)$ is nonsingular. Then*

(a) $\lim_{\varepsilon \to 0+} \{\phi_D((1-\varepsilon)p + \varepsilon\tilde{p}) - \phi_D(p)\}/\varepsilon = \Sigma_{i=1}^n \tilde{p}_i \{\psi_{Di}(p) - q\}$,

(b) $\lim_{\varepsilon \to 0+} \{\phi_A((1-\varepsilon)p + \varepsilon\tilde{p}) - \phi_A(p)\}/\varepsilon = \Sigma_{i=1}^n \tilde{p}_i [\psi_{Ai}(p) - \text{tr}\{H(p)\}^{-1}]$.

**Lemma 3**. *A design measure p is*

(a) *D-optimal if and only if $H(p)$ is nonsingular and $\psi_{Di}(p) \leq q$, $1 \leq i \leq n$,*

(b) *A-optimal if and only if $H(p)$ is nonsingular and $\psi_{Ai}(p) \leq \text{tr}\{H(p)\}^{-1}$, $1 \leq i \leq n$.*

**Lemma 4**. *Let $\delta$ be a preassigned positive quantity, $\hat{p}$ be a design measure with nonsingular $H(\hat{p})$ and $\bar{\phi}_D = \max \phi_D(p)$, $\bar{\phi}_A = \max \phi_A(p)$, the maxima being over all design measures.*

(a) *If $\psi_{Di}(\hat{p}) - q \leq \delta$, $1 \leq i \leq n$, then $\phi_D(\hat{p}) \geq \bar{\phi}_D - \delta$.*

(b) *If $\psi_{Ai}(\hat{p}) - \text{tr}\{H(\hat{p})\}^{-1} \leq \delta$, $1 \leq i \leq n$, then $\phi_A(\hat{p}) \geq \bar{\phi}_A - \delta$.*

Proofs of Lemmas 1 and 2 appear in the appendix. The arguments for Lemma 2 are similar to but more elaborate than those in Mukerjee and Huda (2014) who considered only the *A*-criterion in a different context, effectively with $t = 1$ in $H(p)$. Lemmas 3 and 4 follow from Lemmas 1 and 2 arguing as in Silvey (1980, pp. 19-20, 35).

**Example 1**. Suppose the design points nonnull vectors with elements 0 and $\pm 1$, as happens for instance in chemical balance weighing designs. Let $w$ be the smallest integer such that $w \geq q+1$ and a Hadamard matrix of order $w$ exists; e.g., if $q = 6$, then $w = 8$. Without loss of generality, let the first row of this Hadamard matrix consist only of +1's. Choose any $q$ of the remaining rows to form a $q \times w$ matrix $D$ with elements $\pm 1$ and consider the design measure $\bar{p}$ which assigns a mass $1/w$ to each design point represented by a column of $D$. In view of the properties of a Hadamard matrix, then by (1), $g(\bar{p})$ is a null vector and $H(\bar{p}) = G(\bar{p}) = I_q$, the identity matrix of order $q$. So, by (4) and (5), $\psi_{Di}(\bar{p}) = \psi_{Ai}(\bar{p}) = x_i^T x_i \leq q \ [= \text{tr}\{H(\bar{p})\}^{-1}]$ for every $i$, and the *D*- and *A*-optimality of the design measure $\bar{p}$, irrespective of the value of $t$, is immediate from Lemma 3. □

The optimality of $\bar{p}$ in Example 1 can as well be verified from first principles. In Section 3, we consider an application where Lemmas 3 and 4 are really useful in a far more nontrivial manner.

## 3. Binary design points

Let the design points be nonnull binary vectors of order $q \times 1$, as happens for example in spring balance weighing designs. Then the design space $\Omega$ has cardinality $n = 2^q - 1$. Depending on convenience, we will denote the design points interchangeably by either $x_1, \ldots, x_n$ or $x$, $x \in \Omega$. With this design space, Huda and Mukerjee (1988) investigated *D*- and *A*-optimal design measures under OLS estimation. They worked with $G(p)$, which amounts to taking $t = 0$. In contrast, our design problem concerns $0 \le t < 1$ and this warrants significant additional work as the proofs, all reported in the appendix, reveal. These proofs make use of Lemma 3 but involve tricky manipulation of the terms $\psi_{Di}(p)$ and $\psi_{Ai}(p)$ in order to establish the inequalities there.

### 3.1 Analytical results

We begin by considering the case $q = 2$, where $\Omega$ consists of only three design points. Here the *D*-optimality result looks nice but, quite counterintuitively, the *A*-optimal design measure turns out to be rather involved, thus calling for a separate treatment of this case. The following lemma will be helpful.

**Lemma 5**. *For $0 \le t < 1$, let*

$$u_t(\xi) = 1 - 2t\xi - (3 - 2t)\xi^2 + 4t\xi^3 - 2t^2\xi^4.$$

*Then the equation $u_t(\xi) = 0$ has a unique root in $[1/2, 1)$. Moreover, this unique root, say $\xi_t$, satisfies $1/2 < \xi_t < 2^{-1/2}$.*

**Theorem 1**. *Let $q = 2$ and denote the design points by $x_1 = (1,0)^T$, $x_2 = (0,1)^T$ and $x_3 = (1,1)^T$. Then*

(a) *the design measure $p_D$ assigning a mass $1/3$ to each $x_1, x_2, x_3$ is D-optimal for every $0 \le t < 1$;*

(b) *the design measure $p_A$ assigning a mass $1 - \xi_t$ to each of $x_1, x_2$ and a mass $2\xi_t - 1$ to $x_3$, where $\xi_t$ is as defined in Lemma 5, is A-optimal.*

The sharp contrast between parts (a) and (b) of Theorem 1 is worth noting. For $q = 2$, the same design measure is seen in part (a) to be *D*-optimal for every *t*, while part (b) shows that the *A*-optimal design measure is very much dependent on *t* via $\xi_t$. Although an analytical expression for $\xi_t$ is not available, it can be obtained numerically and is shown in Table 1 for various *t*.

Table 1. *Values of $\xi_t$ for various t*

| $t$ | 0 | 0.1 | 0.2 | 0.3 | 0.4 | 0.5 | 0.6 | 0.7 | 0.8 | 0.9 |
|---|---|---|---|---|---|---|---|---|---|---|
| $\xi_t$ | 0.5774 | 0.5858 | 0.5950 | 0.6051 | 0.6162 | 0.6285 | 0.6423 | 0.6580 | 0.6758 | 0.6948 |

Turning to general $q$, let $\Omega_j (\subset \Omega)$ be the set of design points with $j$ ones, $1 \le j \le q$. Each $\Omega_j$ has cardinality $n_j = \binom{q}{j}$. Consider the following design measures none of which involves *t*.

(i) For even $q = 2m$, $m \geq 2$: (a) $p_{ev1}$ is the design measure assigning a mass $(n_m + n_{m+1})^{-1}$ to each design point in $\Omega_m \cup \Omega_{m+1}$ and a mass 0 to every other design point; (b) $p_{ev2}$ is the design measure assigning a mass $n_m^{-1}$ to each design point in $\Omega_m$ and a mass 0 to every other design point.

(ii) For odd $q = 2m+1$, $m \geq 1$: $p_{odd}$ is the design measure assigning a mass $n_{m+1}^{-1}$ to each design point in $\Omega_{m+1}$ and a mass 0 to every other design point.

Theorems 2 and 3 below present analytical optimality results on these design measures.

**Theorem 2**. *Let $q = 2m$, $m \geq 2$, be even. Then*

(a) *the design measure $p_{ev1}$ is D-optimal if and only if $0 \leq t \leq t_1(q)$, where $t_1(q) = (q+1)/(q+2)$;*

(b) *the design measure $p_{ev2}$ is A-optimal if and only if $0 \leq t \leq t_2(q)$, where*

$$t_2(q) = 1 - \tfrac{1}{2}(q-1)^{-2}[q + \{4(q-1)^2 + q^2\}^{1/2}].$$

**Theorem 3**. *Let $q = 2m+1$, $m \geq 1$, be odd, and $t_0(q) = q/(q+1)$. Then*

(a) *the design measure $p_{odd}$ is D-optimal if and only if $0 \leq t \leq t_0(q)$;*

(b) *the design measure $p_{odd}$ is A-optimal if and only if $0 \leq t \leq t_0(q)$.*

It is interesting to note that for odd $q$, the same design measure turns out to be both D- and A-optimal in Theorem 3 over identical ranges of $t$. The quantities $t_0(q)$, $t_1(q)$ and $t_2(q)$ in Theorems 2 and 3 tend to 1 as $q \to \infty$. Indeed, even for moderate $q$, Theorems 2 and 3 yield design measures, free from $t$, which are optimal over a reasonably wide range of $t$, e.g., $t_1(q) = 0.875$ and $t_2(q) = 0.647$ for $q = 6$, while $t_0(q) = 0.833$ for $q = 5$. This is practically useful because in reality the true value of $t$ is rarely known.

*3.2 Algorithmic results*

For $t$ beyond the ranges covered by Theorems 2 and 3, it is hard to find optimal design measures in a closed form. However, multiplicative algorithms (cf. Zhang and Mukerjee, 2013) can be employed conveniently for their numerical determination. The D-optimality algorithm starts with the uniform measure $p^{(0)} = (1/n,...,1/n)^T$ which assigns a mass $1/n$ to each design point and, in view of Lemma 3(a), finds $p^{(h)} = (p_1^{(h)},..., p_n^{(h)})^T$, $h = 1,2,...$ recursively as

$$p_i^{(h)} = p_i^{(h-1)} \psi_{Di}(p^{(h-1)})/q, \quad 1 \leq i \leq n,$$

till a design measure $p^{(h)}$, satisfying $\psi_{Di}(p^{(h)}) - q \leq \delta$, $1 \leq i \leq n$, is obtained, where $\delta (> 0)$ is a preassigned small quantity. The A-optimality algorithm also starts with the uniform measure and, as suggested by Lemma 3(b), finds $p^{(h)} = (p_1^{(h)},..., p_n^{(h)})^T$, $h = 1,2,$ ... recursively as

$$p_i^{(h)} = p_i^{(h-1)} \psi_{Ai}(p^{(h-1)}) / \text{tr}\{H(p^{(h-1)})\}^{-1}, \quad 1 \leq i \leq n,$$

till a design measure $p^{(h)}$, satisfying $\psi_{Ai}(p^{(h)}) - \text{tr}\{H(p^{(h)})\}^{-1} \leq \delta$, $1 \leq i \leq n$, is obtained, for some preassigned small $\delta$ (> 0). In the stopping rules for both algorithms, we take $\delta = 10^{-10}$. Then by Lemma 4, the terminal design measures are optimal, under the respective criteria, with accuracy up to 9 places of decimals. Even at this level of accuracy, the algorithms are seen to be quite fast for smaller $q$. They can, however, be slow for larger $q$, but this is not of much concern, for then Theorems 2 and 3 cover a very wide range of $t$ and hence there is no serious need for algorithmic determination of optimal design measures.

For even $q$, $4 \leq q \leq 10$, Tables 2 and 3 show D- and A-optimal design measures for $t$ not covered by Theorem 2, i.e., for $t > t_1(q)$ and $t > t_2(q)$, respectively. We consider $t$ up to 0.9, in interval of 0.1. Since $t_1(4) = 0.833$, $t_1(6) = 0.875$ and $t_1(q) \geq 0.9$ for $q \geq 8$, only two cases need to be considered in Table 2: $q = 4, 6$ and $t = 0.9$. For each such $q$ and $t$, Table 2 also shows the D-efficiency of the design measure $p_{\text{ev1}}$, computed as $\text{eff}_D(p_{\text{ev1}}) = [\det\{H(p_{\text{ev1}})\}/\det\{H(\bar{p}_D)\}]^{1/q}$, where $\bar{p}_D$ is the D-optimal design measure. This helps in assessing the performance of $p_{\text{ev1}}$, seen in Theorem 2(a) to be D-optimal over an appreciable range of $t$, when $t$ falls outside this range. For the same purpose, we report in Table 3 the A-efficiency of the design measure $p_{\text{ev2}}$ in Theorem 2(b). This is computed as $\text{eff}_A(p_{\text{ev2}}) = \text{tr}\{H(\bar{p}_A)\}^{-1}/\text{tr}\{H(p_{\text{ev2}})\}^{-1}$, where $\bar{p}_A$ is the A-optimal design measure. Table 3 also shows $t_2(q)$ because its form is a bit involved. Since the convergence of $t_2(q)$ to 1 is not as fast as that of $t_1(q)$, Table 3 is more elaborate than Table 2.

For odd $q$, Table 4 shows D- and A-optimal design measures for $t$ not covered by Theorem 3, i.e., for $t > t_0(q)$. These D- and A-optimal design measures turn out to be identical. Thus the pattern seen in Theorem 3 persists although the design measure $p_{\text{odd}}$ there no longer remains optimal. Again, we consider $t$ up to 0.9, in interval of 0.1. Since $t_0(3) = 0.75$, $t_0(5) = 0.833$, $t_0(7) = 0.875$ and $t_0(q) \geq 0.9$ for $q \geq 9$, only four cases need to be considered in Table 4: $q = 3$, $t = 0.8$, and $q = 3, 5, 7$, $t = 0.9$. Table 4 also exhibits the D- and A-efficiencies of $p_{\text{odd}}$, which are computed as before. In Tables 2-4, $\pi_j$ denotes the mass assigned by the optimal design measure to each design point in $\Omega_j$. Throughout Tables 3 and 4, it is found that $\pi_j = 0$, for $j \geq 7$. Hence in these two tables we show $\pi_j$ only for $j \leq 6$.

The D- and A-efficiency figures in Tables 2-4 are quite impressive, especially for $q \geq 5$ when they are always over 0.90 and often close to 1. Thus, for $q \geq 5$, the design measures in Theorems 2 and 3 are not only optimal over a wide range of $t$ but also have high efficiency even for $t$ outside this range. So, they can be safely used in the absence of explicit knowledge of $t$, which is typically the case in practice.

Table 2. *D-optimal design measures for even q and $t_1(q) < t \leq 0.9$*

| $q$ | $t$ | $\text{eff}_D(p_{\text{ev1}})$ | $\pi_1$ | $\pi_2$ | $\pi_3$ | $\pi_4$ | $\pi_5$ | $\pi_6$ |
|---|---|---|---|---|---|---|---|---|
| 4 | 0.9 | 0.9807 | 0.0444 | 0.0778 | 0.0778 | 0.0444 | | |
| 6 | 0.9 | 0.9968 | 0.0006 | 0.0073 | 0.0241 | 0.0241 | 0.0073 | 0.0006 |

Table 3. *A-optimal design measures for even q and $t_2(q) < t \leq 0.9$*

| $q$ | $t_2(q)$ | $t$ | $\text{eff}_A(p_{\text{ev2}})$ | $\pi_1$ | $\pi_2$ | $\pi_3$ | $\pi_4$ | $\pi_5$ | $\pi_6$ |
|---|---|---|---|---|---|---|---|---|---|
| 4 | 0.377 | 0.4 | 0.9999 | 0.0000 | 0.1644 | 0.0034 | 0.0000 | | |
| 4 | | 0.5 | 0.9974 | 0.0000 | 0.1535 | 0.0197 | 0.0000 | | |
| 4 | | 0.6 | 0.9891 | 0.0000 | 0.1407 | 0.0390 | 0.0000 | | |
| 4 | | 0.7 | 0.9685 | 0.0000 | 0.1253 | 0.0620 | 0.0000 | | |
| 4 | | 0.8 | 0.9190 | 0.0000 | 0.1068 | 0.0898 | 0.0000 | | |
| 4 | | 0.9 | 0.7579 | 0.0445 | 0.0777 | 0.0779 | 0.0443 | | |
| 6 | 0.647 | 0.7 | 0.9991 | 0.0000 | 0.0000 | 0.0459 | 0.0054 | 0.0000 | 0.0000 |
| 6 | | 0.8 | 0.9875 | 0.0000 | 0.0000 | 0.0367 | 0.0177 | 0.0000 | 0.0000 |
| 6 | | 0.9 | 0.9257 | 0.0006 | 0.0074 | 0.0240 | 0.0241 | 0.0073 | 0.0005 |
| 8 | 0.754 | 0.8 | 0.9989 | 0.0000 | 0.0000 | 0.0000 | 0.0125 | 0.0022 | 0.0000 |
| 8 | | 0.9 | 0.9763 | 0.0000 | 0.0000 | 0.0000 | 0.0079 | 0.0079 | 0.0000 |
| 10 | 0.811 | 0.9 | 0.9916 | 0.0000 | 0.0000 | 0.0000 | 0.0000 | 0.0025 | 0.0018 |

Table 4. *D- and A-optimal design measures for odd q and $t_0(q) < t \leq 0.9$*

| $q$ | $t$ | $\text{eff}_D(p_{\text{odd}})$ | $\text{eff}_A(p_{\text{odd}})$ | $\pi_1$ | $\pi_2$ | $\pi_3$ | $\pi_4$ | $\pi_5$ | $\pi_6$ |
|---|---|---|---|---|---|---|---|---|---|
| 3 | 0.8 | 0.9902 | 0.9846 | 0.0625 | 0.2500 | 0.0625 | | | |
| 3 | 0.9 | 0.8842 | 0.8000 | 0.1667 | 0.1111 | 0.1667 | | | |
| 5 | 0.9 | 0.9751 | 0.9529 | 0.0082 | 0.0313 | 0.0481 | 0.0313 | 0.0082 | |
| 7 | 0.9 | 0.9963 | 0.9931 | 0 | 0.0005 | 0.0070 | 0.0170 | 0.0070 | 0.0005 |

*3.3 Reducing the support size*

As seen in the last two subsections, the design measures in Theorems 2 and 3 are either optimal or highly efficient over a very wide range of $t$, particularly when $q \geq 5$. However, their support size increases rapidly with $q$ and, with a finite number of observations, this may hinder their practical implementation. In the spirit of Huda and Mukerjee (1988), we now indicate the use of balanced incomplete block (BIB) designs to find design measures which are exactly or approximately as good as the ones in Theorems 2 and 3 but have much smaller support size. As noted earlier, while Huda and Mukerjee (1988) worked with $G(p)$, i.e., $t = 0$, we consider $H(p)$ for general $t$.

For ease in reference, recall that a BIB$(q, b, r, k, \lambda)$ design, say $d$, is an arrangement of $q$ symbols in $b$ blocks of size $k$ ($< q$) such that every symbol appears at most once in each block, every symbol appears in $r$ blocks, and every two distinct symbols occur together in $\lambda$ blocks. The incidence matrix of $d$ is denoted by $N_d$, which is $q \times b$ with $(i, j)$th entry 1 if the $i$th symbol appears in the $j$th block, and 0 otherwise. Clearly, the $b$ columns of $N_d$ are $q \times 1$ binary vectors, i.e., design points in our design space $\Omega$. Let $p(d)$ be a design measure that assigns a mass $1/b$ to each column of $N_d$. Then from (1), using the standard relationships $qr = bk$ and $r(k-1) = \lambda(q-1)$ for a BIB design, one can check that $H(p(d))$ is uniquely determined by the parameters $q$ and $k$ of $d$.

Now observe that the design measures $p_{ev2}$ and $p_{odd}$ in Theorems 2 and 3 are supported on design points which, when written as columns, form incidence matrices of BIB designs, say $d_{ev2}$ and $d_{odd}$, having parameters $q = 2m, b = \binom{2m}{m}, k = m$ and $q = 2m+1, b = \binom{2m+1}{m+1}, k = m+1$, respectively. Moreover, both $p_{ev2}$ and $p_{odd}$ spread the total mass uniformly on their support points. Hence, as noted above, if we can find a BIB design $d$ which has the same $q$ and $k$ as $d_{ev2}$ or $d_{odd}$ but smaller $b$, then $p(d)$ will have the same $H(.)$ as $p_{ev2}$ or $p_{odd}$ but a smaller support size. From this perspective, we consider the following BIB designs:

$d_1$: $q = 2m, b = 4m-2, r = 2m-1, k = m, \lambda = m-1$;

$d_2$: $q = 4s+1, b = 8s+2, r = 4s+2, k = 2s+1, \lambda = 2s+1$;

$d_3$: $q = b = 4s+3, r = k = 2s+2, \lambda = s+1$.

It is well-known (Raghavarao, 1971, Ch. 5) that $d_1$ is coexistent with a Hadamard matrix of order $4m$, $m \geq 2$, $d_3$ is coexistent with a Hadamard matrix of order $4s+4$, and $d_2$ exists provided $4s+1$ is a prime or prime power. Thus for each $3 \leq q \leq 20$, one of the these BIB designs exist. Moreover, $d_1$ has the same $q$ and $k$ as $d_{ev2}$ and either $d_2$ or $d_3$ has the same $q$ and $k$ as $d_{odd}$ according as $q = 1$ or 3 (mod 4). Hence writing $p^{[j]} = p(d_j)$, $j = 1, 2, 3$, the following are evident:

(i) for even $q = 2m$, $m \geq 2$, $H(p^{[1]}) = H(p_{ev2})$;

(ii) for odd $q = 2m+1$, $H(p^{[2]}) = H(p_{odd})$ if $m = 2s$ is even, and $H(p^{[3]}) = H(p_{odd})$ if $m = 2s+1$ is odd.

By (i), $p^{[1]}$ enjoys the optimality or high efficiency properties of $p_{ev2}$ while entailing significant reduction in support size. For example, with $q$ = 6, 8 and 10, $p^{[1]}$ is supported on 10, 14 and 18 points as against 20, 70 and 252 support points in $p_{ev2}$. Similarly, by (ii) $p^{[2]}$ or $p^{[3]}$ enjoy all the nice properties $p_{odd}$ noted earlier while having a much smaller support size. For instance, with $q$ = 7 and 9, $p^{[3]}$ and $p^{[2]}$ have 7 and 18 support points, respectively, while $p_{odd}$ is supported on 35 and 126 points. Note also that none of $p^{[1]}, p^{[2]}, p^{[3]}$ involves $t$.

The aforesaid technique of reducing the support size does not work for the *D*-optimal design measure $p_{ev1}$ in Theorem 2 because its support points do not form the incidence matrix of a BIB design. However, for $q = 4$, $p_{ev1}$ itself is supported on only 10 points, while for even $q \geq 6$, we find that $p^{[1]}$ comes quite close to $p_{ev1}$ under the *D*-criterion. This is evident from Table 5 which shows the *D*-efficiency of $p^{[1]}$ relative to $p_{ev1}$, as given by $[\det\{H(p^{[1]})\}/\det\{H(p_{ev1})\}]^{1/q}$, for $t = 0(0.1)0.9$. At the same time, $p^{[1]}$ is supported on much fewer points than $p_{ev1}$; e.g., with $q = 6, 8$ and 10, $p^{[1]}$ has only 10, 14 and 18 support points as against 35, 126 and 462 in $p_{ev2}$.

To summarize, for $q \geq 6$, the following points emerge on the basis of our findings on optimality or high efficiency vis-à-vis support size, taking due cognizance of the possible uncertainty about *t*:

(a) If *q* is even, then use $p^{[1]}$ under both *D*- and *A*-criteria. It is equivalent to $p_{ev2}$ under the *A*-criterion and competes very well with $p_{ev1}$ under the *D*-criterion, while having much smaller support size than both.

(b) If *q* is odd, then use $p^{[2]}$ or $p^{[3]}$ under both *D*- and *A*-criteria according as $q = 1$ or 3 (mod 4). They are equivalent to $p_{odd}$ under both criteria but are supported on much fewer design points.

Table 5. *D-efficiency of* $p^{[1]}$ *relative to* $p_{ev1}$

| | | | | | t | | | | | |
|---|---|---|---|---|---|---|---|---|---|---|
| q | 0 | 0.1 | 0.2 | 0.3 | 0.4 | 0.5 | 0.6 | 0.7 | 0.8 | 0.9 |
| 6 | 0.9927 | 0.9924 | 0.9919 | 0.9913 | 0.9905 | 0.9894 | 0.9878 | 0.9851 | 0.9798 | 0.9652 |
| 8 | 0.9968 | 0.9966 | 0.9964 | 0.9961 | 0.9957 | 0.9952 | 0.9945 | 0.9932 | 0.9908 | 0.9837 |
| 10 | 0.9983 | 0.9982 | 0.9981 | 0.9979 | 0.9977 | 0.9975 | 0.9971 | 0.9964 | 0.9950 | 0.9911 |

**Appendix: Proofs**

*Proof of Lemma 1.* By (1), for design measures $p$ and $\tilde{p}$ and $0 < \varepsilon < 1$,

$$G((1-\varepsilon)p + \varepsilon\tilde{p}) = (1-\varepsilon)\{H(p) + tg(p)g(p)^T\} + \varepsilon\{H(\tilde{p}) + tg(\tilde{p})g(\tilde{p})^T\},$$

and $g((1-\varepsilon)p + \varepsilon\tilde{p}) = (1-\varepsilon)g(p) + \varepsilon g(\tilde{p})$. Hence on simplification,

$$H((1-\varepsilon)p + \varepsilon\tilde{p}) = (1-\varepsilon)H(p) + \varepsilon H(\tilde{p}) + t\varepsilon(1-\varepsilon)g_0 g_0^T, \qquad (A.1)$$

where $g_0 = g(\tilde{p}) - g(p)$. Thus $H((1-\varepsilon)p + \varepsilon\tilde{p}) - \{(1-\varepsilon)H(p) + \varepsilon H(\tilde{p})\}$ is nonnegative definite. The lemma now follows from (2) and (3) by the usual arguments in approximate theory; cf. Silvey, 1980, p. 17. □

*Proof of Lemma 2.* Write $e_i = x_i - g(p)$, $1 \leq i \leq n$. With $g_0$ as in (A.1), then by (1),

$$G(\tilde{p}) - g(\tilde{p})g(\tilde{p})^T = \Sigma_{i=1}^n \tilde{p}_i \{x_i - g(\tilde{p})\}\{x_i - g(\tilde{p})\}^T = \Sigma_{i=1}^n \tilde{p}_i (e_i - g_0)(e_i - g_0)^T.$$

Hence for $0 < \varepsilon < 1$, by (1) and (A.1), $H((1-\varepsilon)p + \varepsilon\tilde{p}) = (1-\varepsilon)H(p) + \varepsilon Q(\varepsilon)$, where

$$Q(\varepsilon) = H(\widetilde{p}) + t(1-\varepsilon)g_0 g_0^T = (1-t)G(\widetilde{p}) + t\{G(\widetilde{p}) - g(\widetilde{p})g(\widetilde{p})^T + (1-\varepsilon)g_0 g_0^T\}$$

$$= (1-t)\Sigma_{i=1}^n \widetilde{p}_i x_i x_i^T + t\{\Sigma_{i=1}^n \widetilde{p}_i (e_i - g_0)(e_i - g_0)^T + (1-\varepsilon)g_0 g_0^T\} = U(\varepsilon)U(\varepsilon)^T$$

where $U(\varepsilon)$ is a $q \times (2n+1)$ matrix with columns $\{(1-t)\widetilde{p}_i\}^{1/2} x_i$, $(t\widetilde{p}_i)^{1/2}(e_i - g_0)$, $1 \leq i \leq n$, and $(1-\varepsilon)^{1/2} g_0$. Hence the determinant and inverse of $H((1-\varepsilon)p + \varepsilon \widetilde{p})$ equal

$$\det\{(1-\varepsilon)H(p)\} \det[I_{2n+1} + \tfrac{\varepsilon}{1-\varepsilon} U(\varepsilon)^T \{H(p)\}^{-1} U(\varepsilon)],$$

and

$$\tfrac{1}{1-\varepsilon}[\{H(p)\}^{-1} - \tfrac{\varepsilon}{1-\varepsilon}\{H(p)\}^{-1} U(\varepsilon)\{I_{2n+1} + \tfrac{\varepsilon}{1-\varepsilon} U(\varepsilon)^T (H(p))^{-1} U(\varepsilon)\}^{-1} U(\varepsilon)^T \{H(p)\}^{-1}],$$

respectively. Noting that $Q(\varepsilon)$ and $U(\varepsilon)$ are well-defined also at $\varepsilon = 0$, it now follows from (2) and (3) that

$$\phi_D((1-\varepsilon)p + \varepsilon \widetilde{p}) = \phi_D(p) - \varepsilon q + \varepsilon \text{tr}[U(0)^T \{H(p)\}^{-1} U(0)] + O(\varepsilon^2),$$

$$\phi_A((1-\varepsilon)p + \varepsilon \widetilde{p}) = \phi_A(p) - \varepsilon \text{tr}\{H(p)\}^{-1} + \varepsilon \text{tr}[\{H(p)\}^{-1} U(0) U(0)^T \{H(p)\}^{-1}] + O(\varepsilon^2).$$

Therefore, as $U(0)U(0)^T = Q(0)$, the limits in parts (a) and (b) equal

$$\text{tr}[\{H(p)\}^{-1} Q(0)] - q \quad \text{and} \quad \text{tr}[\{H(p)\}^{-1} Q(0)\{H(p)\}^{-1}] - \text{tr}\{H(p)\}^{-1},$$

respectively. Since $\Sigma_{i=1}^n \widetilde{p}_i e_i = g_0$, and hence $Q(0) = (1-t)\Sigma_{i=1}^n \widetilde{p}_i x_i x_i^T + t\Sigma_{i=1}^n \widetilde{p}_i e_i e_i^T$, the truth of the lemma is now evident from (4) and (5). □

**Proof of Lemma 5.** Since $u_t(1/2) = \tfrac{1}{8}(2 - t^2) > 0$, $u_t(2^{-1/2}) = -\tfrac{1}{2}(1-t)^2 < 0$, and

$$u_t(\xi) = \tfrac{1}{2}\xi^{-2}\{(1-\xi)^2(1-2\xi^2) - (2t\xi^3 - 2\xi^2 - \xi + 1)^2\}, \qquad \xi \in [1/2, 1),$$

we conclude that $u_t(\xi) = 0$ has a root in $(1/2, 2^{-1/2})$ and no root in $[2^{-1/2}, 1)$. It now suffices to show that $u_t(\xi)$ is strictly decreasing in $\xi$ over $[1/2, 2^{-1/2}]$. This follows because for any such $\xi$,

$$u'_t(\xi) = -2t - 2(3 - 2t)\xi + 12t\xi^2 - 8t^2\xi^3$$

$$= -\tfrac{1}{8}\xi^{-3}\{12\xi^2(1-\xi)^2 - (1-2\xi)^2 + (8t\xi^3 - 6\xi^2 - 2\xi + 1)^2\} < 0.$$

Here we use the facts that the minimum of $12\xi^2(1-\xi)^2 - (1-2\xi)^2$, over $\xi$ in $[1/2, 2^{-1/2}]$, is attained at $\xi = 2^{-1/2}$, and that this minimum is positive. □

**Proof of Theorem 1.** (a) By (1), $g(p_D)$ has each element 2/3. Also, upon finding $H(p_D)$ from (1), one can check that it is nonsingular and that its inverse has each diagonal element $3(6-4t)/(9-8t)$ and each off-diagonal element $-3(3-4t)/(9-8t)$. Hence (4) yields $\psi_{Di}(p_D) = 2 (= q)$, for each $i$ and every $0 \leq t < 1$, and the D-optimality of $p_D$ follows from Lemma 3(a).

(b) By (1), $g(p_A)$ has each element $\xi_t$. Also, upon finding $H(p_A)$ from (1), one can check that it is nonsingular and that its inverse has each diagonal element $\xi_t(1-t\xi_t)/\Delta$ and each off-diagonal element $-(2\xi_t - 1 - t\xi_t^2)/\Delta$, where $\Delta = (1-\xi_t)(3\xi_t - 1 - 2t\xi_t^2)$. Hence, after some algebraic manipulation, (5) yields

$$\text{tr}\{H(p_A)\}^{-1} - \psi_{A1}(p_A) = \text{tr}\{H(p_A)\}^{-1} - \psi_{A2}(p_A) = (2\xi_t - 1)u_t(\xi_t)/\Delta^2,$$

$$\text{tr}\{H(p_A)\}^{-1} - \psi_{A3}(p_A) = 2(\xi_t - 1)u_t(\xi_t)/\Delta^2,$$

where $u_t(.)$ is as in Lemma 5. Since $u_t(\xi_t) = 0$ by the definition of $\xi_t$, the A-optimality of $p_A$ follows from Lemma 3(b). □

Two more lemmas are needed for proving Theorems 2 and 3. In what follows, $1_q$ is the $q \times 1$ vector of ones and $J_q = 1_q 1_q^T$.

**Lemma A.1.** (i) *For even* $q = 2m$, $m \geq 2$, $g(p_{\text{ev1}}) = \{(m+1)/(2m+1)\}1_q$ *and* $g(p_{\text{ev2}}) = (1/2)1_q$. *Also, both $H(p_{\text{ev1}})$ and $H(p_{\text{ev2}})$ are nonsingular, and*

$$\{H(p_{\text{ev1}})\}^{-1} = \frac{2(2m+1)}{m+1}\{I_q - q^{-1}J_q + \frac{2m+1}{1+4m(m+1)(1-t)}(q^{-1}J_q)\},$$

$$\{H(p_{\text{ev2}})\}^{-1} = \frac{2(2m-1)}{m}(I_q - q^{-1}J_q) + \frac{2}{m(1-t)}(q^{-1}J_q).$$

(ii) *For odd* $q = 2m+1$, $m \geq 1$, $g(p_{\text{odd}}) = \{(m+1)/(2m+1)\}1_q$. *Also, $H(p_{\text{odd}})$ is nonsingular and*

$$\{H(p_{\text{odd}})\}^{-1} = \frac{2(2m+1)}{m+1}\{I_q - q^{-1}J_q + \frac{1}{2(m+1)(1-t)}(q^{-1}J_q)\}.$$

*Proof.* We indicate the proof only for even $q$ and the design measure $p_{\text{ev1}}$. The proofs for the other cases are similar and simpler. Let $\Sigma_j$ denote sum over design points $x$ which belong to $\Omega_j$. Then one can check that

$$\Sigma_j x = (jn_j/q)1_q, \quad \Sigma_j xx^T = [jn_j/\{q(q-1)\}]\{(q-j)I_q + (j-1)J_q\}, \quad 1 \leq j \leq q. \quad (A.2)$$

For even $q = 2m$, by (1), (A.2) and the definition of $p_{\text{ev1}}$,

$$g(p_{\text{ev1}}) = (n_m + n_{m+1})^{-1}(\Sigma_m x + \Sigma_{m+1} x) = c_0 1_q,$$

$$G(p_{\text{ev1}}) = (n_m + n_{m+1})^{-1}(\Sigma_m xx^T + \Sigma_{m+1} xx^T) = c_1 I_q + c_2 J_q,$$

where

$$c_0 = \{q(n_m + n_{m+1})\}^{-1}\{mn_m + (m+1)n_{m+1}\},$$

$$c_1 = \{q(q-1)(n_m + n_{m+1})\}^{-1}\{m(q-m)n_m + (m+1)(q-m-1)n_{m+1}\},$$

$$c_2 = \{q(q-1)(n_m + n_{m+1})\}^{-1}\{m(m-1)n_m + (m+1)mn_{m+1}\}.$$

Since $q = 2m$, on simplification, $c_0 = 2c_1 = 2c_2 = (m+1)/(2m+1)$. Hence $g(p_{ev1})$ is as stated. Next, the assertion regarding $H(p_{ev1})$ follows after some additional steps, using (1). □

**Lemma A.2.** *For even $q = 2m$, $m \geq 2$, let*

$$l(j,t) = \{(q-1)^2(1-t)^2 - 1\}(j-m)^2 - q(1-t)(j-m),$$

*where $0 \leq t < 1$. Then $l(j,t) \geq 0$, $1 \leq j \leq q$ if and only if $0 \leq t \leq t_2(q)$.*

*Proof.* Note that $l(m,t) = 0$ and $l(m+c,t) < l(m-c,t)$, for every $0 \leq t < 1$ and every positive integer $c$. Hence it suffices to show that $l(j,t) \geq 0$, $m+1 \leq j \leq q$ if and only if $0 \leq t \leq t_2(q)$. Since

$$l(j,t) = (q-1)^2(j-m)^2\{1 - t - \tfrac{1}{2}(j-m)^{-1}(q-1)^{-2}q\}^2 - \{(j-m)^2 + \tfrac{1}{4}(q-1)^{-2}q^2\}$$

and $0 \leq t < 1$, it follows that $l(j,t) \geq 0$ for any fixed $j$ ($m+1 \leq j \leq q$) if and only if

$$t \leq 1 - \tfrac{1}{2}(q-1)^{-2}[(j-m)^{-1}q + \{4(q-1)^2 + (j-m)^{-2}q^2\}^{1/2}].$$

The right-hand side of this inequality is minimum, over $m+1 \leq j \leq q$, at $j = m+1$ and this minimum value equals $t_2(q)$. Hence $l(j,t) \geq 0$, $m+1 \leq j \leq q$ if and only if $0 \leq t \leq t_2(q)$. □

***Proof of Theorem 2.*** Note that for any design point $x \in \Omega_j$ and any constant $c$,

$$(x - c1_q)^T(I_q - q^{-1}J_q)(x - c1_q) = j(q-j)/q,$$

$$(x - c1_q)^T(q^{-1}J_q)(x - c1_q) = (j - cq)^2/q, \tag{A.3}$$

In what follows, we will often use the fact that here $q = 2m$.

(a) Write $c = (m+1)/(2m+1)$. From (4), (A.3) and Lemma A.1(i), if a design point $x_i$ belongs to $\Omega_j$, then $\psi_{Di}(p_{ev1}) = f_1(j)$, where

$$f_1(j) = \frac{2(2m+1)}{(m+1)q}\left[j(q-j) + \frac{2m+1}{1 + 4m(m+1)(1-t)}\{(1-t)j^2 + t(j-cq)^2\}\right].$$

Therefore, by Lemma 3(a), $p_{ev1}$ is $D$-optimal if and only if $f_1(j) \leq q$, $1 \leq j \leq q$. After considerable algebra, one can express $f_1(j)$ as

$$f_1(j) = q - \frac{2(2m+1)\{q+1-(q+2)t\}}{(m+1)\{1 + 4m(m+1)(1-t)\}}(j-m)(j-m-1)$$

Part (a) now follows because $(j-m)(j-m-1) \geq 0$ for $1 \leq j \leq q$, the inequality being strict unless $j = m$ or $m+1$.

(b) By Lemma A.1(i),

$$\operatorname{tr}\{H(p_{ev2})\}^{-1} = 2m^{-1}\{(2m-1)^2 + (1-t)^{-1}\}.$$

Also, as $I_q - q^{-1}J_q$ and $q^{-1}J_q$ are both idempotent, by (5), (A.3) and Lemma A.1(i), if a design point $x_i$ belongs to $\Omega_j$, then $\psi_{Ai}(p_{ev2}) = f_2(j)$, where

$$f_2(j) = 4(m^2 q)^{-1}[(2m-1)^2 j(q-j) + (1-t)^{-2}\{(1-t)j^2 + t(j-m)^2\}].$$

Therefore, by Lemma 3(b), $p_{ev2}$ is A-optimal if and only if $f_2(j) \leq \text{tr}\{H(p_{ev2})\}^{-1}$, $1 \leq j \leq q$. An intricate algebra shows that

$$\text{tr}\{H(p_{ev2})\}^{-1} - f_2(j) = 16 q^{-3}(1-t)^{-2} l(j,t),$$

with $l(j,t)$ as in Lemma A.2. Part (b) now follows from Lemma A.2. □

*Proof of Theorem 3*. Here $q = 2m+1$, a fact which will often be used in the proof.

(a) From (4), (A.3) and Lemma A.1(ii), if a design point $x_i$ belongs to $\Omega_j$, then $\psi_{Di}(p_{odd}) = f_3(j)$, where

$$f_3(j) = \frac{2}{m+1}\left(j(q-j) + \frac{(1-t)j^2 + t(j-m-1)^2}{2(m+1)(1-t)}\right) = q - \frac{\{q-(q+1)t\}(j-m-1)^2}{(m+1)^2(1-t)},$$

upon simplification. Hence part (a) follows from Lemma 3(a).

(b) By Lemma A.1(ii),

$$\text{tr}\{H(p_{odd})\}^{-1} = \frac{2q}{m+1}\left(q - 1 + \frac{1}{2(m+1)(1-t)}\right).$$

Also, as $I_q - q^{-1}J_q$ and $q^{-1}J_q$ are idempotent, by (5), (A.3) and Lemma A.1(ii), if a design point $x_i$ belongs to $\Omega_j$, then $\psi_{Ai}(p_{odd}) = f_4(j)$, where

$$f_4(j) = \frac{4q}{(m+1)^2}\left(j(q-j) + \frac{(1-t)j^2 + t(j-m-1)^2}{4(m+1)^2(1-t)^2}\right).$$

A long and tricky algebra yields

$$\text{tr}\{H(p_{odd})\}^{-1} - f_4(j) = \frac{q\{q-(q+1)t\}\{(j-m-1)^2 + 2(j-m-1)(j-m)(m+1)(1-t)\}}{(m+1)^4(1-t)^2}.$$

Part (b) now follows from Lemma 3(b), since $(j-m-1)^2 + 2(j-m-1)(j-m)(m+1)(1-t) \geq 0$ for $1 \leq j \leq q$, the inequality being strict unless $j = m+1$. □

**Acknowledgement**: The work of RM was supported by the J.C. Bose National Fellowship of the Government of India and a grant from Indian Institute of Management Calcutta.